%% file: Airy_Pearcey_PDE3.tex
\documentclass[12pt]{article}
\usepackage{amsmath}
\usepackage{amssymb}
\usepackage{bbm}
\numberwithin{equation}{section}
\usepackage{epic,pdfsync}

\usepackage{amssymb,color,graphicx}

  \centerline{\LARGE From the Pearcey to the Airy process}
\centerline{\LARGE  }

\author{ }

\bigskip
\centerline{Mark Adler\footnote{Department of Mathematics, Brandeis University, Waltham, Mass 02454, USA,
adler@brandeis.edu. The support of a National Science Foundation grant \# DMS-07-04271 is gratefully acknowledged},
Mattia Cafasso\footnote{Centre de recherches math\'ematiques,
Universit\'e de Montr\'eal; C.~P.~6128, succ. centre ville, Montr\'eal,
Qu\'ebec, Canada H3C 3J7\\ Department of Mathematics and
Statistics, Concordia University; 1455 de Maisonneuve W., Montr\'eal, Qu\'ebec,
Canada H3G 1M8, cafasso@crm.umontreal.ca}, Pierre van Moerbeke\footnote{D\'epartement de Math\'ematiques, Universit\'e Catholique
de Louvain, 1348 Louvain-la-Neuve, Belgium and Brandeis University, Waltham, Mass 02454, USA,
pierre.vanmoerbeke@uclouvain.be. The support of a National Science Foundation grant \# DMS-07-04271, a European Science
Foundation grant (MISGAM), a Marie Curie Grant (ENIGMA), Nato, FNRS and Francqui Foundation grants is gratefully
acknowledged.} }

\date{}

\newcommand{\MAT}[1]{\left(\begin{array}{*#1c}}
\newcommand{\mat}{\end{array}\right)}
\newcommand{\qed}{\leavevmode\unskip\nobreak\penalty200\hskip2pt\null
\nobreak\hfill\rule{1.1ex}{1.1ex}
\medbreak }

\newcommand{\AR}{{\cal A}}

\newcommand{\PR}{{\cal P}}

\newcommand{\BP}{{\mathbb P}}

\newcommand{\iy}{\infty}
\newcommand{\pl}{\partial}

\newenvironment{proof}{\medskip\noindent{\it Proof:\/} }{\qed}

\newcommand{\om}{\omega}

\newcommand{\Dt}{\Delta}
 \newcommand{\vr}{\varepsilon}

\newcommand{\BR}{{\mathbb R}}

\newcommand{\lb}{\lambda}

\newcommand{\dis}{\displaystyle}

 \newcommand{\un}{\mbox{1\hspace{-3.2pt}I}}

\def\be#1\ee{\begin{equation}#1\end{equation}}
\def\bea#1\eea{\begin{eqnarray}#1\end{eqnarray}}
\def\bean#1\eean{\begin{eqnarray*}#1\end{eqnarray*}}

\newcommand{\Tr}{\operatorname{\rm Tr}}

\newtheorem{definition}{Definition}[section]
\newtheorem{theorem}[definition]{Theorem}
\newtheorem{lemma}[definition]{Lemma}

\newtheorem{proposition}[definition]{Proposition}

\catcode `!=11

\newdimen\squaresize
\newdimen\thickness
\newdimen\Thickness
\newdimen\ll! \newdimen \uu! \newdimen\dd! \newdimen \rr! \newdimen
\temp!

\def\sq!#1#2#3#4#5{%
\ll!=#1 \uu!=#2 \dd!=#3 \rr!=#4
\setbox0=\hbox{%
 \temp!=\squaresize\advance\temp! by .5\uu!
 \rlap{\kern -.5\ll!
 \vbox{\hrule height \temp! width#1 depth .5\dd!}}%
 \temp!=\squaresize\advance\temp! by -.5\uu!
 \rlap{\raise\temp!
 \vbox{\hrule height #2 width \squaresize}}%
 \rlap{\raise -.5\dd!
 \vbox{\hrule height #3 width \squaresize}}%
 \temp!=\squaresize\advance\temp! by .5\uu!
 \rlap{\kern \squaresize \kern-.5\rr!
 \vbox{\hrule height \temp! width#4 depth .5\dd!}}%
 \rlap{\kern .5\squaresize\raise .5\squaresize
 \vbox to 0pt{\vss\hbox to 0pt{\hss $#5$\hss}\vss}}%
}
 \ht0=0pt \dp0=0pt \box0
}

\def\vsq!#1#2#3#4#5\endvsq!{\vbox to \squaresize{\hrule
width\squaresize height 0pt%
\vss\sq!{#1}{#2}{#3}{#4}{#5}}}

\newdimen \LL! \newdimen \UU! \newdimen \DD! \newdimen \RR!

\def\vvsq!{\futurelet\next\vvvsq!}
\def\vvvsq!{\relax
  \ifx     \next l\LL!=\Thickness \let\continue=\skipnexttoken!
  \else\ifx\next u\UU!=\Thickness \let\continue=\skipnexttoken!
  \else\ifx\next d\DD!=\Thickness \let\continue=\skipnexttoken!
  \else\ifx\next r\RR!=\Thickness \let\continue=\skipnexttoken!
  \else\def\continue{\vsq!\LL!\UU!\DD!\RR!}%
  \fi\fi\fi\fi
  \continue}

\def\skipnexttoken!#1{\vvsq!}

\def\place#1#2#3{\vbox to 0pt{\vss
\rlap{\kern#1\squaresize
  \raise#2\squaresize\hbox{$#3$}}
\vss}}

\catcode `!=12

\begin{document}

\begin{abstract}

Putting dynamics into random matrix models leads to finitely many nonintersecting Brownian motions on the real line for the eigenvalues, as was discovered by Dyson. Applying scaling limits to the random matrix models, combined with Dyson's dynamics, then leads to interesting, infinite-dimensional diffusions for the eigenvalues. This paper studies the relationship between two of the models, namely the Airy and Pearcey processes and more precisely shows how to approximate the multi-time statistics for the Pearcey process by the one of the Airy process with the help of a PDE governing the gap probabilities for the Pearcey process. 

\end{abstract}


\section{Introduction} 

Putting dynamics into random matrix models leads to finitely many nonintersecting Brownian motions on $\BR$ for the eigenvalues, as was discovered by Dyson \cite{Dy:BrownianMotions}. Applying scaling limits to the random matrix models, combined with Dyson's dynamics, then leads to interesting, infinite-dimensional diffusions for the eigenvalues. This paper studies the relationship between two of the models, namely the Airy and Pearcey processes and more precisely shows how to approximate the Pearcey process by the Airy process with the help of a PDE \cite{AOvM} governing the gap probabilities for the Pearcey process. The Airy process was introduced by Pr\"ahofer-Spohn \cite{PrSp} and further developed by K. Johansson \cite{JoDetProc,Johansson3}. A simple non-linear 3rd order PDE for the transition probabilities for this process was found in \cite{AvM-Airy-Sine}; see also \cite{TW-Airy,TW-Dyson}. 
%
The Pearcey process was introduced in \cite{TracyWidomPearcey,OkounkovReshetikhin} in the context of non-intersecting Brownian motions and plane partitions, also based on prior work on matrix models with external source \cite{MooreMatrix,BowickBrezin,Zinn1,Zinn2,Brezin2,  Brezin5,     AptBleKui,BleKui1,BleKui2,BleherKuijlaarsIII}.


\newpage

Consider $n$ nonintersecting Brownian particles on the real line $\BR$,
$$
-\iy <x_1(t)< ... < x_n(t)<\iy
$$
with (local) Brownian transition probability given by
$$
p(t;x,y):=\frac{1}{\sqrt{\pi t}}e^{-\frac{(y-x)^2}{t}},
$$
all starting from the origin $x=0$ at time $t=0$ and such that $\frac{n}{2}$ particles are forced to $\pm\sqrt{\frac{n}{2}}$ at $t=1$. 
%
%
For very large $n$, the average mean density of particles has its support, for $t\leq\frac{1}{2}$, on one interval centered about $x=0$, and for $ \frac{1}{2}<t\leq 1$ on two intervals symmetrically located about $x=0$. The end points of the interval(s) of support describe a heart-shaped region in $(x,t)$-space, with a cusp at $(x_0,t_0)=(0,1/2)$. The Pearcey process ${\cal P}(\tau)$ is defined (see Figure 1) as the motion of these nonintersecting Brownian motions for large $n$, about $(x_0,t_0)=(0,1/2)$ (i.e., near the cusp), with space microscopically rescaled by a factor of $n^{-1/4}$ and time rescaled by a factor $n^{-1/2}$, in tune with the Brownian motion rescaling. A partial differential equation for the Pearcey process was found by Adler-van Moerbeke \cite{AvM-Pearcey} and in \cite{AOvM} a much better version was obtained, namely a simple third order non-linear PDE for the transition probabilities; it was obtained by a scaling limit on a PDE for non-intersecting Brownian motions with target points. This PDE is related to the Boussinesq equation and its hierarchy; this is part of a general result on integrable kernels, as explained in the forthcoming paper \cite{ACvM}.  

Near the boundary of the heart-shaped region of Figure 1, but away from the cusp, the local fluctuations behave as the so-called Airy process, which describe the non-intersecting brownian motions with space stretched by the customary GUE edge rescaling $n^{1/6}$ and time rescaled by the factor $n^{1/3}$, again in tune with the Brownian motion space-time rescaling. 

This paper shows how the Pearcey process statistics tends to the Airy process statistics when one is moving out of the cusp $x=\frac{2}{27}(3(t-t_0))^{3/2}$ very near the boundary, that is at a distance of $(3\tau)^{1/6}$ for $\tau$ very large, with $\tau$ being the Pearcey time. To be precise, in the two-time case, the times $\tau_i$ must be sufficiently near -in a very precise way- for the limit to hold. 
The main result of the paper can be summarized as follows:

 \begin{theorem}\label{Th: Airy by Pearcey}
 Given finite parameters $t_1<t_2$, let both $\tau_1,~\tau_2 \to \iy$, such that $\tau_2-\tau_1\to \iy$ behaves in the following precise way: 
  \be
  \frac{\tau_2-\tau_1}{2(t_2-t_1)}=(3\tau_1)^{1/3}
  +\frac{t_2-t_1}{(3\tau_1)^{1/3}}+\frac{2t_1t_2}{3\tau_1}+O\bigl( {\tau_1^{-5/3}}\bigr)
  .\label{est-t}
  \ee
 The parameters $t_1$ and $t_2$ provide the Airy times in the following approximation of the Airy process by the Pearcey process:
 \be\begin{aligned}
 \lefteqn{ \BP\left(\bigcap_{i=1}^2\left\{\frac{\PR(\tau_i)-\frac{2}{27}(3\tau_i)^{3/2}}{(3\tau_i)^{1/6}}\cap\left(-E _i
\right)=\emptyset   \right\} \right)} 
 \\
&~~~~~~= 
\BP\left( \bigcap_{i=1}^2\left\{{\cal A}(t_i)\cap (-  E_i)=\emptyset 
\right\}\right) \left(1+O\bigl(\frac{1}{\tau_1^{4/3}}\bigr)\right). 
\end{aligned} \label{est}\ee

  \end{theorem}
  
  The same estimate holds as well for the one-time case. A similar (but different) result was then obtained in \cite{BC} for the one--time case in the situation when one is moving out of the cusp following both the two branches of the cusp simultaneously.

  The proof of this Theorem proceeds in two steps. At first, Proposition \ref{Th: prop 3.5}, stated later, shows that, using the scaling above, the Pearcey kernel tends to the Airy kernel. In a second step, we show, using the PDE for the Pearcey process \cite{AOvM} and Proposition \ref{Th: prop 3.5}, the result of Theorem \ref{Th: Airy by Pearcey}.



It is a fact that both the Airy and Pearcey processes are determinantal processes, for which the multi-time gap probabilities is given by the matrix Fredholm determinant of the matrix kernel, which will be described here.

The {\bf Airy process} is a determinantal process, for which the multi-time gap probabilities\footnote{$\chi_E$ is the indicator function for the set $E$} are given by the matrix Fredholm determinant for $t_1<\ldots <t_\ell$,
\be
\BP\left(\bigcap^{\ell}_{i=1}(\AR(t_i)\cap E_i=\emptyset)\right)=\det\left(I-[\chi_{E_i}K^{\AR}_{t_i t_j}\chi_{E_j}]\right)_{1\leq i,j\leq\ell}\\
\ee
of the matrix kernel in $\vec x=(x_1,\ldots,x_\ell)$ and $\vec y=( y_1,\ldots,y_\ell)$, denoted as follows:
\be\begin{aligned}
\lefteqn{{\mathbb K}^{\AR}_{t_1,\ldots,t_\ell}(\vec x, \vec y)\sqrt{d\vec x~ d\vec y}}
\\ \\
&:=\left(\begin{array}{llll}
  K^{\cal A}_{t_1,t_1}(x_1,y_1)\sqrt{dx_1dy_1} & \ldots & K^{\cal A}_{t_1,t_\ell} (x_1,y_\ell) \sqrt{dx_1dy_\ell}\\ 
  \vdots &  & \vdots \\
  K^{\cal A}_{t_\ell,t_1} (x_\ell,y_1)\sqrt{dx_\ell dy_1}& \ldots & K^{\cal A}_{t_\ell,t_\ell} (x_\ell,y_\ell)\sqrt{dx_\ell dy_\ell}
  \end{array}\right) ,
\end{aligned}\ee
where, for arbitrary $t_i$ and $t_j$, the extended Airy kernel $K^{\AR}_{t_i,t_j}(x,y) $ is given by (see Johansson \cite{JoDetProc,Johansson3})
\be
K^{\AR}_{t_i,t_j}(x,y):=\tilde K^{\AR}_{t_i,t_j}(x,y)-   \un  (t_i<t_j)p^{\cal A} (t_j-t_i,x,y),
\label{extAiry kernel}\ee
  where 
  $$
   \begin{aligned}
  {\tilde K^{\AR}_{t_i,t_j}(x,y) } 
 &:=  \int_0^{\iy} e^{-\lb(t_i-t_j)}{  A}(x+\lb){  A}(y+\lb)d\lb      \\
&~=   \frac{1}{\left(2\pi i\right)^2}
\int_{\Gamma_{>}}\!
dv  \int_{\Gamma_{<}}\! du ~  \frac{e^{-v^3/3+yv}}{e^{- u^{3}/3+x  u}}
\frac{1}{(v+t_j)-( u+t_i)}
\\ \\
p^{\cal A} (t,x,y)&:=  \frac{1}{\sqrt{4\pi t}} e^{\frac{t ^{3}}{12}-
   \frac{\left( x-y \right) ^{2}}{4 t}-\frac{t}{2}   \left( y+x \right) } ,\mbox{ for } t>0.
\end{aligned} 
$$
The contour $u\in \Gamma_{<}$ consists of two rays emanating from the origin with angles $\theta_1$ and $\theta'_1$ 
  with the positive real axis, 
  and the contour $v \in \Gamma_{>}$ also consists of two rays with angles $\theta_2$ and $\theta'_2$ 
   with the negative real axis, as indicated in Figure 2. As is well known, one may choose $\theta_1, \theta_2,~\theta'_1,\theta'_2=  \pi/3 $.

In particular, for $t_i=t_j$, this is the customary Airy kernel
\be
K^{\cal A}_{t_i,t_i}:=K^{\cal A}(x,y):=\int_0^\iy d\lb { A}(x+\lb)
{  A}(y+\lb) =\frac{A(x)A'(y)-A(y)A'(x)}{x-y}.
\label{Airy kernel}
\ee
%


\begin{figure}[center]
\hspace{1.6cm}
\resizebox{0.8\textwidth}{!}{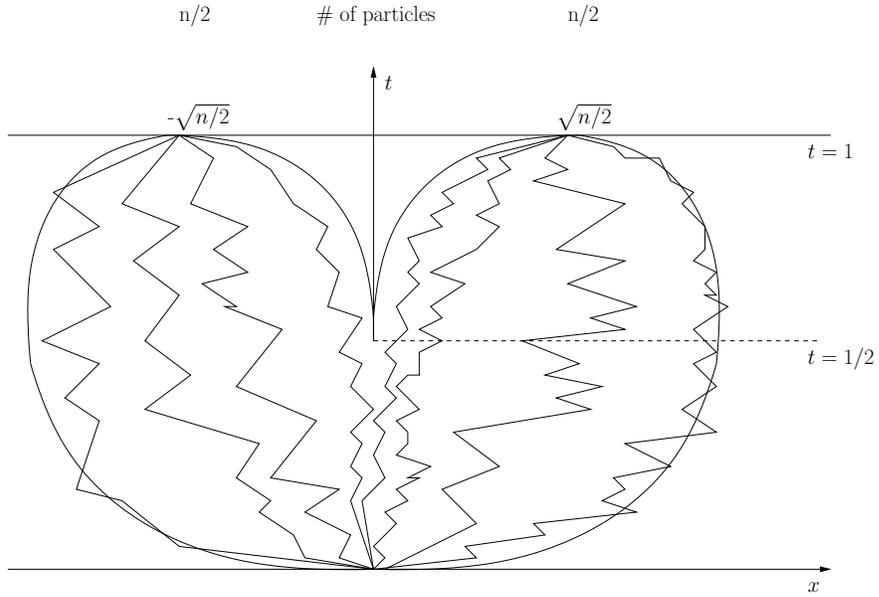}
\caption{The Pearcey process.}
\end{figure}

  
 \vspace*{5cm}
 
  $ \hspace*{-4cm}   \begin{picture}(-12,1.6)
\put(220,70){\circle*{1}}

\multiput(220,70)(0,10){8}{\line(0,1){8}}
  \multiput(220,70)(0,-10){8}{\line(0,-1){8}}

\put(180,110){\vector(-1,1){3}}
\put(180,30){\vector( 1, 1){3}} 

\put(260,110){\vector( -1,-1){3}}
\put(260,30){\vector( 1, -1){3}} 

\put(213,110){\vector(  -1,4){1}}
\put(213,30){\vector( 1, 4){1}}

 \put(220,70){\line(-1,1){50}} 
 \put(220,70){\line(-1,-1){50}}
    \put(220,70){\line(-1,6){14}}
    \put(220,70){\line(-1,-6){14}}
 
  
    \multiput(220,70)(10, 0){8}{\line(1,0){8}}
    \multiput(220,70)(-10, 0){8}{\line(-1,0){8}}
 
   \put(170,30){\makebox(0,0)[r]{$U\in X_2$}}
     \put(170,110){\makebox(0,0)[r]{$U\in X_2$}}
    
     \put(190,50){\makebox(0,0)[r]{$ \sigma'_2$}}
 \put(190,90){\makebox(0,0)[r]{$\sigma_2$}}

   \put(215,60){\makebox(0,0)[r]{$ \tau'$}}
 \put(215,80){\makebox(0,0)[r]{$\tau $}}

 \put(270,30){\makebox(0,0)[l]{$ U\in X_1$}}
 \put(270,110){\makebox(0,0)[l]{$ U\in X_1$}}
 
  \put(205,-10){\makebox(0,0)[r]{$V\in Y$}}
     \put(205,150){\makebox(0,0)[r]{$V\in Y$}}

  \put(250,50){\makebox(0,0)[l]{$ \sigma'_1$}}
 \put(250,90){\makebox(0,0)[l]{  $ \sigma_1 $}}

 \put(220,70){\line( 1,1){50}}
  \put(220,70){\line( 1,-1){50}}
  
 \put(160,-40){\makebox(0,0)[l]{ ~~     \footnotesize  $\pi/8<\sigma_1,\sigma_2, \sigma'_1,\sigma'_2< 3\pi/8 $ }}
 

\put(160,-60){\makebox(0,0)[l] {~~~\footnotesize $ 3\pi/8<\tau,\tau'< 5\pi/8$}}
  
  \put(180,-90){\makebox(0,0)[l]{Pearcey contour}}
  

 \end{picture}
 $ 
 $  \hspace*{7cm}   
 \begin{picture}(-12,1.6)
\put(220,70){\circle*{1}}

 \put(199,110){\vector(-1,3){3}}
\put(201,30){\vector( 1, 3){3}} 

\put(239,110){\vector( -1,-3){3}}
\put(241,30){\vector( 1, -3){3}} 


 \put(220,70){\line(-1,2){30}} 
 \put(220,70){\line(-1,-2){30}}
 
 \multiput(220,70)(0,10){8}{\line(0,1){8}}
  \multiput(220,70)(0,-10){8}{\line(0,-1){8}}
  
   \multiput(220,70)(10, 0){8}{\line(1,0){8}}
   \multiput(220,70)(-10, 0){8}{\line(-1,0){8}}
 
   \put(190,20){\makebox(0,0)[r]{$v\in  \Gamma_{>}$}}
    \put(190,120){\makebox(0,0)[r]{$v\in  \Gamma_{>} $}}
    
     \put(190,50){\makebox(0,0)[l]{$ \theta'_2$}}
 \put(190,90){\makebox(0,0)[l]{$ \theta_2$}}

 \put(250,20){\makebox(0,0)[l]{$ u\in  \Gamma_{<} $}}
 \put(250,120){\makebox(0,0)[l]{$ u \in  \Gamma_{<} $}}
 
  \put(240,50){\makebox(0,0)[l]{  $  \theta'_1 $}}
 \put(240,90){\makebox(0,0)[l]{  $ \theta_1 $}}

 \put(220,70){\line( 1,2){30}}
  \put(220,70){\line( 1,-2){30}}

 \put(160,-40){\makebox(0,0)[l]{ ~      \footnotesize  $\pi/6<\theta_1,\theta_2,\theta'_1,\theta'_2<  \pi/2 $ }}
 
  
  \put(180,-90){\makebox(0,0)[l]{Airy contour}}
  
   \put(10,-120){\makebox(0,0)[l]{Figure 2: Integration paths for the kernels.}}

  
 
  \end{picture}
 $ \hspace{15cm} 
%
%

 \vspace*{5cm}
 \newpage



 
 



The {\bf Pearcey process} is determinantal as well, for which the multi-time gap probability is given by the Fredholm determinant  for $\tau_1<\ldots <\tau_\ell$,
\be 
\BP\left(\bigcap^{\ell}_{i=1}(\PR(\tau_i)\cap E_i=\emptyset)\right)=\det\left(\un -[\chi_{E_i}K^{\PR}_{\tau_i\tau_j}\chi_{E_j}]\right)_{1\leq i,j\leq\ell} 
\label{ProbP}\ee
for the {\bf Pearcey matrix kernel} in $\vec \xi=(\xi_1,\ldots,\xi_\ell)$ and $\vec \eta=( \eta_1,\ldots,\eta_\ell)$, denoted as follows:
\be\begin{aligned}
 \lefteqn{{\mathbb K}^{\cal P}_{\tau_1,\ldots,\tau_\ell}(\vec \xi, \vec \eta)\sqrt{d\vec \xi~ d\vec \eta}}\\ \\
 &:=\left(\!\begin{array}{lllll}
  K^{\cal P}_{\tau_1,\tau_1} (\xi_1	,\eta_1)\sqrt{d\xi_1 d\eta_1} & \ldots & K^{\cal P}_{\tau_1,\tau_\ell}(\xi_1,\eta_\ell) \sqrt{d\xi_1 d\eta_\ell}  \\   
  \vdots & & \vdots \\
  K^{\cal P}_{\tau_\ell,\tau_1}(\xi_\ell,\eta_1)\sqrt{d\xi_\ell d\eta_1}    & \ldots & K^{\cal P}_{\tau_\ell,\tau_\ell} (\xi_\ell,\eta_\ell)\sqrt{d\xi_\ell d\eta_\ell} 
  \end{array}\!\right) ,
\end{aligned}\ee
where for arbitrary $\tau_i$ and $\tau_j$, (see Tracy-Widom \cite{TracyWidomPearcey})
\be
  K^{\PR}_{\tau_i,\tau_j}(\xi,\eta)=\tilde K^{\PR}_{\tau_i,\tau_j}(\xi,\eta) -\un (\tau_i<\tau_j)p^{\cal P}(\tau_j-\tau_i;\xi,\eta),   
\label{Pearcey kernel} \ee
with
$$
  \begin{aligned}
\tilde K^{\PR}_{\tau_i,\tau_j}(\xi,\eta)  
&:= 
- \frac{ 1}{4\pi^2 }\int_X dU
\int_Y
 \frac{dV }{V\!-\!U}~~
\frac{e^{-\frac{V^4}{4}+\frac{\tau_jV^2}{2}-V\eta}}
{ e^{-\frac{U^4}{4}+\frac{\tau_iU^2}{2}-U\xi}}  \\
%
 p^{\cal P}(\tau;\xi,\eta) & :=\frac{1}{\sqrt{2\pi\tau}}e^{-\frac{(\xi-\eta)^2}{2\tau}}
,\mbox{ for } \tau>0 ,\end{aligned}
  $$
 where $X$ and $Y$ are the following contours: $X=X_1\cup X_2$ consists of four rays emanating from the origin with angles $\sigma_1, \sigma'_1$ with the positive real axis and $\sigma_2,\sigma'_2$ with the negative real axis, as given in Figure 2. The contour $Y$ consists of two rays emanating from the origin with angles $\tau, \tau'$ with the negative real axis; it is customary to pick $\sigma_1=\sigma_2= \sigma'_1= \sigma'_2=\pi/4$ and $\tau=\tau'=\pi/2$.  


In \cite{AvM-Pearcey,AOvM}, it was shown that, given intervals   $$
 E_i=(\xi^{(i)}_1,\xi^{(i)}_2),
 $$
 the log of the probability
$$
Q=Q(\tau_1,\ldots,\tau_\ell; E_1,\ldots,E_\ell):=\log \BP\left(\bigcap^{\ell}_{i=1}\{\PR(\tau_i)\cap E_i=\emptyset\}\right)
$$ satisfies a non-linear PDE, which we describe here. Given the times and the intervals above, one defines the following operators:
 $$
 \pl_\tau:=\sum_i\frac{\pl}{\pl \tau_i}  ,\qquad  \pl_{_{\!E_i}}:=\sum_k\frac{\pl}{\pl \xi^{(i)}_k} 
   , \qquad\pl_{_{\!E}}:=\sum_i \pl_{_{\!E_i}}$$
 
 \be
  \vr_\tau:=\sum_i \tau_i\frac{\pl}{\pl \tau_i}  
, \qquad \vr_E:=\sum_i  \sum_k \xi^{(i)}_k\frac{\pl}{\pl \xi^{(i)}_k} .
 \label{operators}\ee
 Then $Q$ satisfies the Pearcey partial differential equation in its arguments and the boundary of the intervals:
 \be
 2\pl_\tau^3Q+\frac{1}{4}(2\vr_\tau+\vr_E-2)\pl_E^2Q
 -   (\sum_i \tau_i\pl_{_{E_i}})   \pl_\tau\pl_EQ
 + \left\{ \pl_\tau\pl_E Q,\pl_E^2Q \right\}_{_{\!\pl_E}}=0
 .\label{PDE}\ee

%




\section{From the Pearcey to the Airy kernel }\label{sect2}

 \noindent Define the rational functions 
      \be \begin{aligned}  
\Phi(x,t;u) &:= -\frac{1}{  4(3u)^{3}}-\frac{t}{(3u)^2}+\frac{x-t^2}{3u}+\frac{4}{3} tx+\frac u6 t^2x \\
 h(x,t;u)&:=\frac{ux}{4}(x+6t^2)
 \label{Phi} %
,\end{aligned}  \ee  
and the diagonal matrix 
$$
S:=\left( \begin{array}{ccl}
 e^{\Phi(x,t_1;z^4)-h(x,t_1;z^4) } \!\!& \!\! 0\\
 0\!\!& \!\!  e^{\Phi(x,t_2;z^4)-h(x,t_2;z^4) }
 \end{array}\right).
 $$
 We now state:
 
 
\begin{proposition}\label{Th: prop 3.5}
Given a parameter $z\to 0$, define two times $\tau_1$ and $\tau_2$ blowing up like $z^{-6}$ and depending on two parameters $t_1$ and $t_2$, 
\be \tau_i=\frac{1} 
{3z^6}(1+6t_iz^4) + O(z^{10})\label{T-scaling}
.\ee
Given large $\xi_i$ and $\eta_j$, define new space variables $x_i$ and $y_j$, using $\tau_i$ above,
\be\begin{aligned}
\xi_i&= \frac{2}{27}(3\tau_i)^{3/2}-(3\tau_i)^{1/6}x _i\\
\eta _j&= \frac{2}{27}(3\tau_j)^{3/2}-(3\tau_j)^{1/6}y _j
\label{X-scaling}
.\end{aligned}\ee
With this space-time rescaling, the following asymptotic expansion holds for the Pearcey kernel in powers of $z^4$, with polynomial coefficients in $t_1+t_2$,  
  \be
  S~
  {\mathbb K}^{\cal P}_{\tau_1,\tau_2}(\vec \xi,\vec \eta)
 \sqrt{d\vec \xi d\vec \eta}~ S^{-1} 
 =\left(1
  +(t_1\!+\!t_2)O_1(z^4)\right)
   {\mathbb K}^{\cal A}_{t_1,t_2}(\vec x,\vec y)\sqrt{\vec dx \vec dy} +O_2(z^8),\label{1.15}\ee
 where $O_1$ refers to a differential operator in $\pl /\pl x,\pl/\pl y$, with polynomial coefficients in $x,y,t_1\pm t_2$ and $(t_1-t_2)^{-1}$ and varying with the four matrix entries of the matrix kernel.

\end{proposition}

On general grounds, due to the Fredholm determinant formula, this would give Theorem \ref{Th: Airy by Pearcey}, but with a much poorer estimate in (\ref{est}), namely instead of $O( {\tau_1^{-4/3}})$, the estimate would be $O( {\tau_1^{-2/3}})$. The existence of the PDE for the Fredholm determinant of the Pearcey process enables us to bootstrap the $O( {\tau_1^{-2/3}})$ estimate to $O( {\tau_1^{-4/3}})$, without tears!

 Also note that conjugating a kernel does not change its Fredholm determinant. 
 Before proving Proposition \ref{Th: prop 3.5}, we shall need the following identities:

%
%



\begin{lemma}\label{Lemma 3.5}

Introducing the polynomial
\be
\Psi (x,s;\om) := 
  \frac14  \om^4 +\frac32  \om^2s^2 
 -4s(x\om-\om^3)
 \label{Psi} %
,\ee
the Airy kernel (\ref{Airy kernel}) and the (double) integral part of the extended Airy kernel (\ref{extAiry kernel}), at $t_1=s,~t_2=-s$, satisfy the following differential equations,
\be \begin{aligned}\lefteqn{
\Bigl(\Psi (x, s,-\pl_x)-\Psi(y, s,\pl_y) + 4s\Bigr)K^{\cal A}(x,y)}
 \\
& \qquad\qquad\qquad\qquad\qquad
=\frac14(x-y)(x+y+6s^2 )K^{\cal A}(x,y) ,
\end{aligned} \label{Identity1}\ee 
and
 \be\begin{aligned}
\lefteqn{ \left(\Psi(x,  s,-\pl_x)-\Psi(y,- s,\pl_y)-\frac32 s\frac{\pl}{\pl s}(\pl_x-\pl_y)\right) \tilde K^{\AR}_{_{{s,-s}
}}(x,y)} \\
&  \qquad\qquad\qquad\qquad\qquad
= \frac14(x - y)(x+y+6s^2) \tilde K^{\AR}_{_{{s,-s}
}}(x,y).
 \end{aligned} \label{Identity2}\ee 
 
 \end{lemma}


 \proof The operator on the left hand side of (\ref{Identity1}) reads 
 \be\begin{aligned} \label{2.4}
\lefteqn{ \Psi((x, s,-\pl_x)-\Psi((y,  s,\pl_y) + 4s } \\
&=&
  4s(1+y\pl_y-\pl_y^3+x\pl_x-\pl_x^3)+\frac32 s^2(\pl_x^2-\pl_y^2)+\frac14
(\pl_x^4-\pl_y^4)
.\end{aligned}\ee
 Using the first representation (\ref{Airy kernel}) of the Airy kernel, one checks using the differential equation for the Airy kernel, $A''(x)=xA(x)$ and thus $A'''(x)=xA'(x)+A(x)$ and $A^{(iv)}(x)=2A'(x)+x^2A(x)$, and using differentiation by parts to establish the last equality,
\bean
\lefteqn{\left((y\pl_y-\pl_y^3)+(x\pl_x-\pl_x^3)\right)
 K^{\cal A}(x,y)}\\
 &=& -\int_0^\iy  dz\left(z\frac{d}{dz}+2\right) A(x+z)A(y+z)
 =-K^{\cal A}(x,y)
. \eean
In order to take care of the other pieces in (\ref{2.4}), one uses the second representation (\ref{Airy kernel}) of the kernel $K^{\cal A}(x,y)$, yielding
 \bean
\bigl( \pl_x^2-\pl_y^2\bigr)K^{\cal A}(x,y)&=& (x-y)K^{\cal A}(x,y)\\
\bigl( \pl_x^4-\pl_y^4\bigr)K^{\cal A}(x,y)&=& (x^2-y^2)K^{\cal A}(x,y)
.\eean
This establishes the first identity (\ref{Identity1})
  of Lemma \ref{Lemma 3.5}. The operator on the left hand side of the second identity (\ref{Identity2}) reads:
  \bean
 \lefteqn{ \Psi(x,  s,-\pl_x)-\Psi(y,-s,\pl_y)-\frac32 s\frac{\pl}{\pl s}(\pl_x-\pl_y)
  }\\
  &=&
 \!\!\! \frac{1}{4}(\pl_x^4\!-\!\pl_y^4)+\frac32 s^2(\pl_x^2\!-\!\pl_y^2)
  + 4s(x\pl_x-\pl_x^3-  y\pl_y+\pl_y^3))-\frac32 s\frac{\pl}{\pl s}(\pl_x\!-\!\pl_y)
  ,~~\eean
 of which we will evaluate all the different terms acting on the integral. Notice at first that, since the expression under the differentiation vanishes at $0$ and $\iy$, one has
\be\begin{aligned} 
0&=  \int_0^\iy dz \frac{\pl}{\pl z}\left(ze^{- 2sz}K^{\cal A}(x+z,y+z)\right)
 \\
&= \int_0^\iy dz ~e^{- 2s z}  K^{\cal A} (x+z,y+z)
- \int_0^{\iy}zdz~
   e^{- 2s z}{ A}(x+z){  A}(y+z) 
 \\
&  - 2s
 \int_0^\iy z dz ~e^{- 2s z}  K^{\cal A} (x+z,y+z),
\end{aligned}\label{A1}\ee  
Again, using the second representation (\ref{Airy kernel}) of the kernel $K^{\cal A}(x,y)$, one checks
 $$\begin{aligned}
  \bigl( \pl_x-\pl_y\bigr)\tilde K^{\AR}_{_{{s,-s}
  }}(x,y)
   &=  -(x-y)\int_0^\iy dz ~e^{- 2s z}  K^{\cal A} (x+z,y+z)
  \\
   \bigl( \pl_x^2-\pl_y^2\bigr)\tilde K^{\AR}_{_{{s,-s}
    }}(x,y)
    &=  (x-y)\tilde K^{\AR}_{_{{s,-s}
    }}(x,y)
 ,\end{aligned}$$
 and also
 \bean 
    s \frac{\pl}{\pl s}\bigl( \pl_x\!-\!\pl_y\bigr)\tilde K^{\cal A}_{_{{s,-s}}}(x,y)  
 &=&   2s(x-y)\int_0^\iy z dz ~e^{- 2s z}  K^{\cal A} (x+z,y+z)
\eean
and, using the differential equation $\bigl(x \pl_x-\pl_x^3\bigr)A(x)=-A(x)$ and (\ref{A1}),
 \bean
   \lefteqn{  - 2s\left( \bigl(x \pl_x-\pl_x^3\bigr)-\bigl(y\pl_y-\pl_y^3\bigr) \right)
  \tilde K^{\cal A}_{_{{s,-s}}}(x,y)  }\\ 
  &=&     - 2s (x-y)\int_0^\iy z dz ~e^{- 2s z}  K^{\cal A} (x+z,y+z)
  \eean
  Using  $A^{(iv)}(x)=2A'(x)+x^2A(x)$ and the Darboux-Christoffel representation (\ref{Airy kernel})  of the Airy kernel, one checks
\bean
\lefteqn{\bigl( \pl_x^4-\pl_y^4\bigr)\tilde K^{\cal A}_{_{{s,-s}}}(x,y)
}\\
&=&  
 \int_0^{\iy}
   e^{- 2s z}({ A}^{(iv)}(x+z){  A}(y+z)-A(x+z){ A}^{(iv)}(y+z))dz 
\\
&=&  
 -2(x\!-\!y)\int_0^\iy  dz ~e^{- 2s z}  K^{\cal A} (x+z,y+z)  
    +(x^2-y^2)\tilde K^{\cal A}_{_{{s,-s}}}(x,y)
   \\
   && +2(x-y)\int_0^\iy z dz ~e^{- 2s z}  { A}(x+z){ A}(y+z)
   .\eean
   Adding all these different pieces and using the expression for $$ - 2s (x-y)\int_0^\iy z dz ~e^{- 2s z}  K^{\cal A} (x+z,y+z),
$$ given by (\ref{A1}),
   leads to the statement of Lemma \ref{Lemma 3.5}.\qed
   
   For the sake of notational convenience in the proof below, set   
\be
{\mathbb K}^{\cal P}_{\tau_1,\tau_2}(\vec \xi, \vec \eta)\sqrt{d\vec \xi~ d\vec \eta}=
     %
    \left(\begin{array}{llll}
  K^{\cal P}_{ 11}(\xi_1,\eta_1) \sqrt{d\xi_1d\eta_1}&   K^{\cal P}_{12}(\xi_1,\eta_2) \sqrt{d\xi_1d\eta_2} \\   
    &    \\
  K^{\cal P}_{21}(\xi_2,\eta_1) \sqrt{d\xi_2 d\eta_1} &   K^{\cal P}_{22} (\xi_2,\eta_2)\sqrt{d\xi_2d\eta_2}
  \end{array}\right) 
\label{Pearcey kernel'}\ee
and similarly for the Airy kernel $ {\mathbb K}^{\cal A}_{t_1t_2}$; the $ \tilde K^{\cal P}_{ij}$ and $ \tilde K^{\cal A}_{ij}$ refer, as before, to the (double) integral part.

  \medskip\noindent{\it Proof of Proposition \ref{Th: prop 3.5}:\/}  Notice that acting with $\pl_x$ and $\pl_y$ on the kernels $\widetilde K^{\cal A}_{_{{11}\atop{22}}} (x,y)$ and $\widetilde K^{\cal A}_{_{{12}\atop{21}}} (x,y)$, amounts to multiplication of the integrand of the kernels with $-u$ and $v$ respectively; i.e., $u \leftrightarrow -\pl_x$ and $v \leftrightarrow \pl_y.$

  Given the (small) parameter $z\in \BR$, consider the following $(i,j)$-dependent change of integration variables 
  $(U,V)\mapsto (u,v)$ in the four Pearcey kernels $K^{\cal P}_{ij}$ in (\ref{Pearcey kernel'}), namely
  %
  %
 \be U=\frac{1}{3z^3}(1+3uz^4)(1+3t_iz^4),~~~V=\frac{1}{3z^3}(1+3 vz^4)(1+3t_jz^4),
\label{2.6'}\ee
together with the changes of variables $(\xi,\eta,\tau_i,\tau_j)\mapsto (x,y,t_i,t_j)$, in accordance with (\ref{T-scaling}) and (\ref{X-scaling}),
%
\be \begin{array}{lll} \dis \tau_i  =\frac{1} 
{3z^6}(1+6t_iz^4)+O(z^{10}), ~~ &
\dis \tau_j =\frac{1} 
{3z^6}(1+6t_jz^4)+O(z^{10})
\\  \\  \dis
\xi = \frac{2}{27}(3\tau_i)^{3/2}-(3\tau_i)^{1/6}x,~~ 
&\dis \eta = \frac{2}{27}(3\tau_j)^{3/2}-(3\tau_j)^{1/6}y
.\end{array}
\label{2.6''}\ee
To be precise, for $K^{\cal P}_{kk}$, one sets in the transformations above $i=j=k$, for $K^{\cal P}_{12}$ and $K^{\cal P}_{21}$, one sets $i=1,~j=2$ and $i=2,~j=1$ respectively.
 Then, remembering the expressions 
 $\Phi$ and $\Psi$, defined in (\ref{Phi}) and (\ref{Psi}), the following estimate holds for small $z$:
\be\begin{aligned}
\frac{e^{-\Phi(y,t_j;z^4)}}{e^{-\Phi(x,t_i;z^4)}}\frac{e^{-\frac{V^4}{4}+\frac{T_jV^2}{2}-V\eta}}
{ e^{-\frac{U^4}{4}+\frac{T_iU^2}{2}-U\xi}}
&= \frac{e^{-z^4\Psi(y,t_j; v)}}{e^{-z^4\Psi(x,t_i;u)}}\frac{e^{-\frac{ v^3}3+y v}}{e^{-\frac{ u^{3}}3+x   u}}
 (1\!+\!t_iO(z^8)\!+\!t_jO(z^8)) 
 \\
 &= 
 \frac{e^{-z^4\Psi(y,t_j;\pl_y)}}{e^{-z^4\Psi(x,t_i;-\pl_x)}}\frac{e^{-\frac{ v^3}3+y v}}{e^{-\frac{ u^{3}}3+x  u}}
 (1\!+\!t_iO(z^8)\!+\!t_jO(z^8)) ;
\end{aligned}   \label{4.X}\ee
replacing in the latter expression $u$ and $v$ by differentiations $\pl_x$ and $\pl_y$ has the advantage that upon doubly integrating in $ u$ and $ v$, the fraction of exponentials in front can be taken out of the integration. 
 Moreover, setting $t_1=t+s$ and $t_2=t-s$, one also checks\footnote{The precise expression for $r(x,y,s):= (x - y)^2-8s^2(x + y-2s^2) $ will be irrelevant in the sequel.}
 :
  \be\begin{aligned}
\lefteqn{\frac{e^{-\Phi(y,t_j;z^4)
 }}{e^{-\Phi(x,t_i;z^4)
 }} p^{\cal P}(\tau_j-\tau_i;\xi,\eta)
 \sqrt{d\xi d\eta}
}\\
&=    p^{\cal A}(-2s,x,y)\sqrt{dxdy}
 \left(1 + \frac{ z^4}{4} \big[(x\!-\!y)(x\!+\!y+6s^2)
  +\frac ts  r(x,y,s)\big]+O(z^8)\right)
%
 .\end{aligned}\label{4.Y}\ee
 %
The multiplication by the quotient of the exponentials on the left hand side of the expressions above will amount to a conjugation of the kernel, which will not change the Fredholm determinant.
Setting $t_i=t+s$ and $t_j=t-s$, with $t=0$, one finds for the non-exponential part in the kernel
\bean
\lefteqn{\sqrt{d\xi d\eta}\frac{dUdV}{V-U}}\\
&=&
  \frac{\sqrt{dxdy}~dudv (1+\frac72(t_i+t_j)z^4)}{v+t_j
   -u-t_i+3z^4(vt_j-ut_i)}+O(z^8)
\\ \\
&=&
\left\{\begin{array}{l}
   \frac{\sqrt{dxdy}~dudv }{v 
   -u }\left(1\pm 4sz^4+O(z^8)\right)
    \mbox{~~for $\left\{\begin{array}{ll} i=1,&j=1\\ 
                                                                      i=2,&j=2 \end{array}\right.$}
   \\ \\
    \frac{\sqrt{dxdy}~dudv }{v 
   -u\mp 2s}\left(1\pm \frac{3z^4s(u+v)}{v-u\mp 2s}+O(z^8)\right)\mbox{~~for $\left\{\begin{array}{ll} i=1,&j=2\\ 
                                                                      i=2,&j=1 .\end{array}\right.$}
  \end{array}\right.
 \eean\be \label{4.Z}\ee
 
%
 Along the same vein as the remark in the beginning of the proof of this Theorem, one notices that multiplication of the integrand of the kernel $\tilde K^{\cal A}_{_{{12}\atop{21}}} (x,y)$ with the fraction, appearing in the last formula of (\ref{4.Z}), amounts to an appropriate differentiation, to wit:
 \be\pm \frac{3z^4s(u+v)}{v-u\mp 2s}\longleftrightarrow
  \frac32 s\frac{\pl}{\pl s}(\pl_y-\pl_x). \label{4.Z1} \ee 
 So we have, using (\ref{4.X}), (\ref{4.Y}), (\ref{4.Z}), and the differential equation for $K^{\cal A}$
  of Lemma \ref{Lemma 3.5}, 
 \be\begin{aligned} 
 \lefteqn{\frac{e^{-\Phi(y,t_{1\atop 2};z^4)}}{e^{-\Phi(x,t_{1\atop 2};z^4)}}\left. K^{\cal P}_{_{{11}\atop{22}}} (\xi,\eta)\sqrt{d\xi d\eta}\right|_{t=0}
 -K^{\cal A}(x,y)\sqrt{dxdy}}           
 \\
     &= 
  z^4\left(
     \pm 4s+\Psi((x,\pm s,-\pl_x)-\Psi((y,\pm s,\pl_y)  \right) K^{\cal A}(x,y) \sqrt{dxdy}+O(z^8)
    \\
  &=    \frac{z^4}4  ~  (x-y)(x+y+6s^2)
     K^{\cal A}  (x,y) \sqrt{dxdy}+O(z^8);
\end{aligned}\label{est1}\ee
the upper(lower)-indices correspond to the upper(lower)-signs. Using (\ref{4.Z1}), and the differential equation for $ \tilde K^{\cal A}_{_{{12}\atop{21}}}$ of Lemma \ref{Lemma 3.5},
 \be\begin{aligned}
 \lefteqn{\frac{e^{-\Phi(y,t_{2\atop 1};z^4)}}{e^{-\Phi(x,t_{1\atop 2};z^4)}}\left. \tilde K^{\cal P}_{_{{12}\atop{21}}} (\xi,\eta)\sqrt{d\xi d\eta}
 -\tilde K^{\cal A}_{_{{12}\atop{21}}}(x,y)\sqrt{dxdy}   \right|_{t=0}} 
       \\
    &= 
     z^4\Bigl(\frac32 s\frac{\pl}{\pl s}(\pl_y\!-\!\pl_x)+\Psi(x,\pm s,-\pl_x)-\Psi(y,\mp s,\pl_y) \Bigr) \left.\tilde K^{\cal A}_{_{{12}\atop{21}}} (x,y)\sqrt{dxdy} \right|_{t=0}\\
    & \qquad\qquad\qquad\qquad\qquad\qquad\qquad\qquad\qquad\qquad+O(z^8)
  \\
  &=   \left. \frac{z^4}4  ~  (x-y)(x+y+6s^2) \tilde K^{\cal A}_{_{{12}\atop{21}}}(x,y)\sqrt{dxdy} \right|_{t=0}+O(z^8).
\end{aligned}\label{est2}\ee
Also, using (\ref{4.Y}),
 \be\begin{aligned}
\lefteqn{\frac{e^{-\Phi(y,t_2;z^4)}}{e^{-\Phi(x,t_1;z^4)}} p(\tau_2-\tau_1;\xi,\eta)\sqrt{d\xi d\eta}
-p^{\cal A}(-2s,x,y)\sqrt{dxdy}\Bigr|_{t=0}
}\\
&=    
  \frac{ z^4}{4} (x-y)(x+y+6s^2) p^{\cal A}(-2s,x,y)\sqrt{dxdy}     +O(z^8).
\end{aligned}\label{est3}\ee
Since $h(y,t_i,z^4)-h(x,t_j;z^4)=
 -\frac{z^4}{4}(x-y)(x+y+6s^2)$ for arbitrary $t_{i,j}=t\pm s$ at $t=0$, and thus
 \be
 \frac{e^{-h(y,t_i;z^4)}}{e^{-h(x,t_j;z^4)}}=1+\frac{z^4}{4} (x-y)(x+y+6s^2)+O(z^8)
 \label{2.16}\ee
 Then the following approximations follow upon combining the three estimates above (\ref{est1}), (\ref{est2}), (\ref{est3}) and using (\ref{2.16}):
$$
{\frac{e^{-\Phi(y,t_{1\atop 2};z^4)+h(y,t_{1\atop 2} ;z^4)}}{e^{-\Phi(x,t_{1\atop 2};z^4)+h(x,t_{1\atop 2} ;z^4)}}\left. K^{\cal P}_{_{{11}\atop{22}}} (\xi,\eta)\sqrt{d\xi d\eta}\right|_{t=0}
 -K^{\cal A}(x,y)\sqrt{dxdy}}=O(z^8)           
$$
 $$
 {\frac{e^{-\Phi(y,t_{2\atop 1};z^4)+h(y,t_{2 \atop 1} ;z^4)}}{e^{-\Phi(x,t_{1\atop 2};z^4)+h(x,t_{1\atop 2} ;z^4)}}\left. \tilde K^{\cal P}_{_{{12}\atop{21}}} (\xi,\eta)\sqrt{d\xi d\eta}
  \right|_{t=0} -\tilde K^{\cal A}_{_{{12}\atop{21}}}(x,y)\sqrt{dxdy}  }  = O(z^8)
$$
$$
 {\frac{e^{-\Phi(y,t_2;z^4)+h(y,t_{2 } ;z^4)}}{e^{-\Phi(x,t_1;z^4)+h(x,t_{1 } ;z^4)}} p^{\cal P}(\tau_2\!-\!\tau_1;\xi,\eta)\sqrt{d\xi d\eta}\Bigr|_{t=0}\!
-p^{\cal A}(-2s,x,y)\sqrt{dxdy}
}=   
   O(z^8).
$$
The reader is reminded that this estimate sofar is done at $t=0$. It then follows for $t\neq 0$ from the estimates (\ref{2.6'}), (\ref{2.6''}), (\ref{4.X}), (\ref{4.Y}), (\ref{4.Z}), (\ref{4.Z1}) that the right hand side of (\ref{1.15}) is an asymptotic series in $z^4$ with polynomial coefficients in $t_1+t_2=2t$, proving the claim about $O_1$.

%

So far an important point was omitted, namely to analyze how the Pearcey contour turns in the limit into an Airy contour; a detailed description of the possible contours was given in Figure 2. The Pearcey-rays $X_1$, with angles $\sigma_1,\sigma'_1 \in (\pi/8,3\pi/8)$ can be deformed into two acceptable rays $\theta_1,\theta'_1\in (\pi/6,\pi/2)$ for the $\Gamma_{<}$-Airy contour, since the two intervals have a non-empty intersection. Also, since $(3\pi/8,5\pi/8)\cap  (\pi/6,\pi/2)\neq \emptyset$, the Pearcey $Y$-contour can be deformed into an acceptable $v \in \Gamma_{>}$-Airy contour. Again, since the admissible interval $(\pi/8,3\pi/8)$ for $\sigma_2$ and $\sigma'_2$ in the $X_2$-Pearcey contour contains $\pi/6$, one ends up in the limit integrating the function $e^{u^3/3}$ along a contour of the form $\Gamma_{>}$, which we may choose to have an angle of $\pi/6-\delta$ with the negative real axis, for small arbitrary $\delta>0$. In the sector $\Gamma_{>}$, with $0<\theta_2= \theta'_2\leq \pi/6-\delta$, the function $e^{u^3/3}$ decays exponentially fast. Therefore, the contribution, due to the $u$-integration of $e^{u^3/3}\times($lower order terms$)$ over $\Gamma_{>}$, having an angle of $\pi/6-\delta$ with the negative real axis,  vanishes by applying Cauchy's Theorem in that sector. 
This ends the proof of Proposition \ref{Th: prop 3.5}.\qed

\section{From the Pearcey to the Airy statistics 
}

This section concerns itself with proving Theorem \ref{Th: Airy by Pearcey}.

\medskip\noindent{\it Proof of Theorem \ref{Th: Airy by Pearcey}:\/}
 For any 
  intervals $E_1$ and $E_2$, the function
   \be
 Q(\tau_1,\tau_2; E_1,E_2)
 =\log \BP\left(\bigcap^{2}_{i=1}(\PR(\tau_i)\cap E_i=\emptyset)\right)
\ee
satisfies the Pearcey PDE (\ref{PDE}). Reparametrizing, without loss of generality, 
%
\be
\tau_1=\tau+\sigma,~~\tau_2=\tau-\sigma,~~E_1= (\xi+\eta+\mu,\xi+\eta-\mu), ~~~E_2=(\xi-\eta+\nu,\xi-\eta-\nu)
 ,\label{int}\ee
 leads to a manageable PDE for the function
\be
F(\tau,\sigma;\xi,\eta,\mu,\nu):=
 Q(\tau_1,\tau_2; E_1,E_2), \label{F}
\ee
namely the PDE:
 \be \begin{aligned}
\lefteqn{\hspace*{-.1cm}2\frac{\pl^3F}{\pl \tau^3}
+\frac14\left(2(\sigma \frac{\pl}{\pl \sigma}\!-\!\tau\frac{\pl}{\pl \tau})
+\xi\frac{\pl}{\pl \xi}+\eta\frac{\pl}{\pl \eta}+\mu\frac{\pl  }{\pl \mu}+\nu\frac{\pl  }{\pl \nu}\!-\!2\right)
 \frac{\pl^2F}{\pl \xi^2}}\\
 & \qquad\qquad\qquad\qquad\qquad
  -     \sigma   \frac{\pl^2F}{\pl \tau\pl \xi \pl \eta}   
 +
  \left\{ \frac{\pl^2F}{\pl \tau\pl \xi},   \frac{\pl^2F}{\pl \xi^2}\right\}_{\xi}
=0.
\end{aligned}  \label{Pearcey-simple}\ee 
To go from the Pearcey PDE (\ref{PDE}) to the PDE (\ref{Pearcey-simple}), one notices that the operators (\ref{operators}), appearing in the Pearcey PDE (\ref{PDE}), have simple expressions in terms of the variables $\tau,\sigma,\xi,\eta,\mu,\nu$,
$$
 \pl_\tau Q=\frac{\pl F}{\pl \tau} ,\qquad \pl_EQ=\frac{\pl F}{\pl \xi} ,
  \qquad \vr_E Q= \bigl(\xi\frac{\pl  }{\pl \xi}+\eta\frac{\pl  }{\pl \eta}+\mu\frac{\pl  }{\pl \mu}+\nu\frac{\pl  }{\pl \nu}\bigr)F,
 $$
 $$  \vr_\tau Q:= \bigl(\sigma\frac{\pl}{\pl \sigma}+\tau\frac{\pl}{\pl \tau} \bigr)F ,~~~ \sum_i \tau_i\pl_{_{\!E_i}} Q= \bigl(\tau\frac{\pl  }{\pl \xi}+\sigma\frac{\pl  }{\pl \eta}\bigr)F
 .$$


The change of variables, considered in (\ref{T-scaling}) and (\ref{X-scaling}), combined with a linear change of variables, in parallel with (\ref{int}), 
\be\begin{aligned}
 \tau_i&= \frac{1} 
{3z^6}(1+6t_iz^4)\mbox{~~~with~~} \left\{\begin{array}{l}t_1=t+s\\t_2=t-s\end{array}\right.
\\
E_i&= \frac{2}{27}(3\tau_i)^{3/2}-(3\tau_i)^{1/6}\tilde E_i
 \mbox{~~~with~~} \left\{\begin{array}{l}\tilde E_1=(x+y+u,x+y-u)\\  \tilde E_2=(x-y+v,x-y-v)\end{array}\right.
\end{aligned}\label{trsf} \ee 
yields a $z$-dependent invertible map,
\be T:~(t,s,x,y,u,v)\mapsto (\tau,\sigma,\xi,\eta,\mu,\nu) \label{map},\ee
and thus the function $F$ in (\ref{F}) leads to a new $z$-dependent function $G$:
$$
F(\tau,\sigma;\xi,\eta,\mu,\nu) =F(T(t,s,x,y,u,v))=:G(t,s,x,y,u,v)
,$$
of which we compute, in principle, the series in $z$. From Proposition \ref{Th: prop 3.5}, it follows that the $z^4$-term in the asymptotic expansion in (\ref{1.15}) vanishes when $t=0$; so, omitting $\sqrt{dx dy}$, one has
$$
S{\mathbbm K}^{\cal P}_{\tau_1,\tau_2}S^{-1}={\mathbbm K}^{\cal A}_{t_1,t_2}
+ z^4 {\mathbbm K}_1+\sum_2^\iy z^{4i}{\mathbbm K}_i
,~~\mbox{with}~  {\mathbbm K}_1= t{\mathbbm H}
,$$
for some kernel ${\mathbbm H}$, with the ${\mathbbm K}_i$ polynomial in $t$. Then defining 
\be
 L_i:=  (I-{\mathbbm K}^{\cal A}_{t_1,t_2})^{-1}{\mathbbm K}_i
 ,\ee
one finds:
\be 
\begin{aligned}
\lefteqn{G(t,s,x,y,u,v)}\\
&={F(\tau,\sigma;\xi,\eta,\mu,\nu)}\\
&  ={\log \BP\Bigl(\bigcap_{i=1}^2\{{\cal P}(\tau_i)\cap E_i=\emptyset \}\Bigr)}\\
&= 
 \log \det (I-{\mathbbm K}^{\cal P}_{\tau_1,\tau_2})_{_{  E_1\times  E_2}}\\
  &= \log \det  (I-{\mathbbm K}^{\cal A}_{t_1,t_2})_{_{\tilde E_1\times \tilde E_2}}
  -z^4\Tr L_1-z^8\Tr (L_2+\frac12   L_1^2)-\ldots
  \\
  &=:
  G_0(s;\tilde E_1,\tilde E_2)-G_1(t,s;\tilde E_1,\tilde E_2)z^4-G_2(t,s;\tilde E_1,\tilde E_2)z^8+O(z^{12}), \label{3.19}
\end{aligned} 
 \ee 
 where \be \begin{aligned}
G_0(s;\tilde E_1,\tilde E_2)&= \log \det  (I-{\mathbbm K}^{\cal A}_{t_1,t_2})_{_{\tilde E_1\times \tilde E_2}} ,\\  \\
G_1(t,s;\tilde E_1,\tilde E_2)&=t\Tr ((I- {\mathbbm K}^{\cal A}_{t_1,t_2})^{-1}{\mathbbm H})_{_{\tilde E_1\times \tilde E_2}} .\label{t-dep}
\end{aligned}\ee 
Note $G_0(s;\tilde E_1,\tilde E_2)$ is $t$-independent, since the Airy process is stationary. Then setting 
$$
F(\tau,\sigma;\xi,\eta,\mu,\nu)
= G(T^{-1}(\tau,\sigma;\xi,\eta,\mu,\nu))=G(t,s,x,y,u,v)
$$
into the PDE (\ref{Pearcey-simple}) yields, using differentiation by parts, a new PDE for $
 G(t, s, x, y, u, v)$. Upon setting the expansion (\ref{3.19}) of $G$ for small $z$, 
$$
 G(t, s, x, y, u, v)=G_0(s;\tilde E_1,\tilde E_2)-G_1(t,s;\tilde E_1,\tilde E_2)z^4+O(z^{8}),
 $$
into this new PDE leads to
  \be
z^2\left(2s\frac{\pl^3 G_0}{\pl s\pl x^2}-\frac{\pl^3 G_1}{\pl t\pl x^2} \right)+O(z^6)=0
,\label{expans1}\ee
with $G_1$ being polynomial in $t$. From the term $\frac{\pl^3F}{\pl \tau^3}$ in the PDE (\ref{Pearcey-simple}), it would seem like the leading term would be of order $z^{-6}$; in fact there are two consecutive cancellations, so that the first non-trivial term has order $z^2$, thus leading to a simple PDE connecting $G_1$ with $G_0$. But since $G_0$ is $t$-independent,  equation (\ref{expans1}) can also be written as 
\be
\frac{\pl}{\pl t}\frac{\pl^2}{\pl x^2}\left(G_1-2ts\frac{\pl}{\pl s} G_0\right) =0
\label{B}\ee
From (\ref{t-dep}), one finds 
 $$\begin{aligned}
G_1-2ts\frac{\pl}{\pl s} G_0&=t\left(\Tr ((I- {\mathbbm K}^{\cal A}_{t_1,t_2})^{-1}{\mathbbm H})_{_{\tilde E_1\times \tilde E_2}}-2s\frac{\pl}{\pl s} G_0\right)\\
&=\sum_1^\ell t^i a_i(s;x,y,u,v),
\end{aligned}
 $$
  which substituted back in (\ref{B}), leads to the PDE's for the $a_i$, namely
  \be
  \left(\frac{\pl}{\pl x}\right)^2
   a_i(s;x,y,u,v) =0
   ,\label{D}\ee
   implying 
   \be a_i(s;x,y,u,v)=xb_i(s;y,u,v)
    +c_i(s;y,u,v).
   \label{E}\ee
  Letting the intervals $\tilde E_1$ and $\tilde E_2$ go to $\iy$, while keeping their relative position $2y$ fixed and widths $-2u$ and $-2v$ fixed as well, is achieved by letting $x\to \iy$, as follows from (\ref{int}). Remember $t_1,t_2, \tilde E_1,\tilde E_2$ from (\ref{trsf}). But in the limit $x\to \iy$, the expressions $G_0$, namely
  $$
  G_0(s;\tilde E_1,\tilde E_2) = \log \det  (I-{\mathbbm K}^{\cal A}_{t_1,t_2})_{_{\tilde E_1\times \tilde E_2}}
  $$ tends to $0$ exponentially fast using the exponential decay of the four components of the matrix Airy kernel (\ref{extAiry kernel}); note that the term $p^{\cal A}(t,x,y)$ tends to $0$ as well, when $x$ and $y$ tend to $\iy$, due to the presence of $e^{-t(x+y)/2}$.  Next, we sketch the proof that
  $$
  G_1(t,s;\tilde E_1,\tilde E_2) =t\Tr ((I- {\mathbbm K}^{\cal A}_{t_1,t_2})^{-1}{\mathbbm H})_{_{\tilde E_1\times \tilde E_2}}
  $$
 tends to $0$ exponentially fast, when $x\to \iy$. Indeed, using the identity
 (${\bf R}$ is the resolvent of the Airy kernel ${\mathbbm K}^{\cal A}_{t_1,t_2} $)
 $$
   (I-{\mathbbm K}^{\cal A}_{t_1,t_2})^{-1}{\mathbbm H}={\mathbbm H}+{\mathbbm K}^{\cal A}_{t_1,t_2}(I-{\mathbbm K}^{\cal A}_{t_1,t_2})^{-1} {\mathbbm H} 
  = 
  {\mathbbm H}+{\mathbbm K}^{\cal A}_{t_1,t_2}(I+{\bf R}) {\mathbbm H},
  $$
  one computes
  $$\begin{aligned}
  G_1&=\Tr ((I-{\mathbbm K}^{\cal A}_{t_1,t_2})^{-1}{\mathbbm K}_1)_{_{\tilde E_1\times \tilde E_2}}
  \\
  &=t\Tr ((I-{\mathbbm K}^{\cal A}_{t_1,t_2})^{-1}{\mathbbm H})_{_{\tilde E_1\times \tilde E_2}}
 \\
 &=t\Tr {\mathbbm H}_{_{\tilde E_1\times \tilde E_2}}+t\Tr ( {\mathbbm K}^{\cal A}_{t_1,t_2}\tilde {\mathbbm H})_{_{\tilde E_1\times \tilde E_2}}, 
  \mbox{  with  } \tilde {\mathbbm H} :=(I+{\bf R}) {\mathbbm H},
  \end{aligned}
  $$
 where
  $$
  \begin{aligned}
  \Tr {\mathbbm H}_{_{\tilde E_1\times \tilde E_2}} &=\sum_{i=1}^2\int_{\tilde E_i} H_{ii}(u,u)du
    \\
  \Tr ({\mathbbm K}^{\cal A}_{t_1,t_2}\tilde {\mathbbm H})_{_{\tilde E_1\times \tilde E_2}} &=\sum_{i=1}^2
  \int\!\!\!\!\!\int_{\tilde E_i\times \tilde E_i}K_{t_it_i}^{\cal A}(u,v)\tilde H_{ii}(v,u)dudv
  \\
  &+
   \int\!\!\!\!\!\int_{\tilde E_1\times \tilde E_2}K_{t_1t_2}^{\cal A}(u,v)\tilde H_{21}(v,u)dudv
   \\&+
   \int\!\!\!\!\!\int_{\tilde E_2\times \tilde E_1}K_{t_2t_1}^{\cal A}(u,v)\tilde H_{12}(v,u)dudv.
  \end{aligned}
  $$
Each of these integrals tend to $0$ because each of them contains the Airy kernel and also because $ {\mathbbm H}$ is obtained by acting with differential operators on the Airy kernel, as explained in section \ref{sect2}. This ends the proof that $
  G_1(t,s;\tilde E_1,\tilde E_2) \to 0$, when the $\tilde E_i$ tend to $\iy$. 
  This fact together with the form (\ref{E}) of the $a_i$ imply $b_i(s;y,u,v)
   =c_i(s;y,u,v)=0$ for $i\geq 1$, and thus $a_i=0$ for $i\geq 1$, implying 
  $$
G_1= 2ts\frac{\pl}{\pl s} G_0.
$$ 
  Summarizing, this implies that for $E_i=\frac{2}{27}(3\tau_i)^{3/2}-(3\tau_i)^{1/6} \tilde E_i$, substituting $G_1=2ts ( \pl G_0/\pl s) $ into (\ref{3.19}),
  \bean
 \lefteqn{ \log \BP\left(\bigcap_{i=1}^2\left\{\frac{\PR(\tau_i)-\frac{2}{27}(3\tau_i)^{3/2}}{(3\tau_i)^{1/6}}\cap\left(-\tilde E _i
\right)=\emptyset   \right\} \right)}\\
& =& G_0(s; \tilde E_1,\tilde E_2) -2z^4 ts\frac{\pl G_0}{\pl s}  +O(z^8)
  \\
&  = & G_0(s-2tsz^4; \tilde E_1,\tilde E_2) +O(z^8) 
 \\
 &=&
 \log \BP\left( \bigcap_{i=1}^2\left\{{\cal A}(t_i(1-t_iz^4))\cap (-\tilde E_i)=\emptyset 
\right\}\right)  +O(z^8),
 \eean
  in view of the change of variables given in (\ref{trsf}), one has that $s-2tsz^4=\frac12(t_1(1-t_1z^4)-t_2(1-t_2z^4))$. Then, substituting $t_i=u_i(1+u_i z^4)$ into $t_i(1-t_iz^4)$ and $\tau_i=(1+6t_iz^4)/(3z^6)+O(z^{10})$, as in (\ref{T-scaling}), yields
 $$
 t_i(1-t_iz^4)= u_i-2u_i^3z^8-u_i^4z^{12} \mbox{  and  }
 \tau_i=\frac1{3z^6} +\frac{2u_i}{z^2}+2u_i^{2} z^2+O(z^{10})
 \mbox{  for  } i=1,2.$$ 
  Then eliminating $z$ between the two expressions above for $\tau_1$ and $\tau_2$, by first expressing $z$ as a series in $\tau_1$, yields (\ref{est-t}) and (\ref{est}), with $t_i$ replaced by $u_i$. This ends the proof of Theorem \ref{Th: Airy by Pearcey}. \qed
  
 \bibliographystyle{plain}
\bibliography{/Users/mattiacafasso/Documents/BibDeskLibrary.bib}

\end{document}

%% file: Pearcey.pdf_t
\begin{picture}(0,0)%
\includegraphics{Pearcey.pdf}%
\end{picture}%
\setlength{\unitlength}{3947sp}%
\begingroup\makeatletter\ifx\SetFigFont\undefined%
\gdef\SetFigFont#1#2#3#4#5{%
  \reset@font\fontsize{#1}{#2pt}%
  \fontfamily{#3}\fontseries{#4}\fontshape{#5}%
  \selectfont}%
\fi\endgroup%
\begin{picture}(10824,7834)(1189,-7052)
\put(5251,539){\makebox(0,0)[lb]{\smash{{\SetFigFont{17}{20.4}{\rmdefault}{\mddefault}{\updefault}{\color[rgb]{0,0,0}$\#$ of particles}%
}}}}
\put(11701,-3961){\makebox(0,0)[lb]{\smash{{\SetFigFont{17}{20.4}{\rmdefault}{\mddefault}{\updefault}{\color[rgb]{0,0,0}$t=1/2$}%
}}}}
\put(11701,-1261){\makebox(0,0)[lb]{\smash{{\SetFigFont{17}{20.4}{\rmdefault}{\mddefault}{\updefault}{\color[rgb]{0,0,0}$t=1$}%
}}}}
\put(11701,-6961){\makebox(0,0)[lb]{\smash{{\SetFigFont{17}{20.4}{\rmdefault}{\mddefault}{\updefault}{\color[rgb]{0,0,0}$x$}%
}}}}
\put(6151,-361){\makebox(0,0)[lb]{\smash{{\SetFigFont{17}{20.4}{\rmdefault}{\mddefault}{\updefault}{\color[rgb]{0,0,0}$t$}%
}}}}
\put(3451,539){\makebox(0,0)[lb]{\smash{{\SetFigFont{17}{20.4}{\rmdefault}{\mddefault}{\updefault}{\color[rgb]{0,0,0}n/2}%
}}}}
\put(8551,539){\makebox(0,0)[lb]{\smash{{\SetFigFont{17}{20.4}{\rmdefault}{\mddefault}{\updefault}{\color[rgb]{0,0,0}n/2}%
}}}}
\put(3301,-811){\makebox(0,0)[lb]{\smash{{\SetFigFont{17}{20.4}{\rmdefault}{\mddefault}{\updefault}{\color[rgb]{0,0,0}-$\sqrt{n/2}$}%
}}}}
\put(8401,-811){\makebox(0,0)[lb]{\smash{{\SetFigFont{17}{20.4}{\rmdefault}{\mddefault}{\updefault}{\color[rgb]{0,0,0}$\sqrt{n/2}$}%
}}}}
\end{picture}%